\documentclass[11pt,a4paper]{article}
\usepackage{amsmath}
\usepackage{amssymb}
\usepackage{amscd}
\usepackage[latin1]{inputenc}
\usepackage[english]{babel} 
\usepackage{fancyhdr}
\usepackage{amsmath,amsfonts}

\def\exp{{\rm exp}\,}

\newtheorem{lemma}{Lemma}[section]
\newtheorem{corollary}[lemma]{Corollary}
\newtheorem{theorem}[lemma]{Theorem}
\newtheorem{proposition}[lemma]{Proposition}
\newtheorem{definition}[lemma]{Definition}
\newtheorem{remark}[lemma]{Remark}

\newenvironment{proof}{\bf Proof : \rm}{\hfill
$\blacksquare$\par\medbreak} 
\date{}

\DeclareMathOperator{\ep}{\varepsilon}

\title{Existence and uniqueness of optimal maps on  Alexandrov spaces}
\author{Jérôme Bertrand\thanks{email: bertrand@math.ups-tlse.fr}}
\begin{document}
\maketitle
\abstract{The purpose of this paper is to show that in a finite dimensional metric space with Alexandrov's curvature bounded below, Monge's transport problem for the quadratic cost  admits a unique solution.}

\section{Introduction}

In this paper, we provide a solution to the Monge-Kantorovich problem in Alexandrov space when the cost function is the square of the metric. We begin by explaining what Monge's problem is, and briefly describing Alexandrov spaces. The mass transportation problem, raised by Monge in 1781 \cite{monge}, is to move one distribution of mass onto another as efficiently as possible, where the criterion of efficiency is to minimise a certain cost. In the original formulation, the cost was the average distance covered by the mass. In other words, given a set $X$ endowed with a cost function, and given two probability measures $\mu_0,\mu_1$  on $X$ as constraints, the problem is to determine, if it exists, a minimizer of
$$ \mbox{ inf } \int_X c(x,s(x))\, d\mu (x)$$
among all measurable maps $s$ sending the initial measure $\mu_0$ onto $\mu_1$ (we denote this condition by $\mu_1=s_{\sharp}\mu_0$).

Monge's original problem has only been solved quite recently (Sudakov gave a proof first \cite{sudakov} but later,  part of his proof appeared to be incorrect,  we refer to \cite{AP} for details). In the meantime, the importance of the mass transport problem for other cost functions was recognized and much work has been devoted to its study. One of the most interesting cases is that of the quadratic cost (the square of the metric distance). In Euclidean space, Monge's problem for the quadratic cost was solved in the eighties by Brenier \cite{brenier} (and  Knott-Smith \cite{KS2} independently) motivated by fluid mechanics problems. Brenier showed the existence and the uniqueness of the optimal map, under the hypothesis that the initial measure is absolutely continous with respect to the Lebesgue measure. The optimal map is the gradient of a convex function on $\mathbb{R}^n$. Subsequently, this result was generalized to the setting of Riemannian manifolds  by McCann \cite{mccann1}, taking advantage of a notion of generalized convex function (so-called $c$-concave function, see Definition \ref{cost} below).

A difficulty in solving Monge's problem is that the problem may be ill-posed, a measure preserving map between two distributions of mass may not exist. To avoid this problem,  Kantorovich introduced a relaxed version of Monge's problem, which is to minimize the quantity 
$$ \int_{X \times X} c(x,y) \; d \Pi(x,y)$$
among all plans, namely all probability measures $\Pi$ whose marginals are $\mu_0$ and $\mu_1$ (in this case, the mass is allowed to split). This problem admits solutions in very general settings (see Section 3). A solution $s$ of Monge's problem for $\mu_0$ and $\mu_1$  induces a solution of Kantorovich problem. More precisely, the plan defined by $(Id,s)_{\sharp}\mu_0$ is a solution. Classically, in the quadratic cost case, existence and uniqueness of the  optimal map follow from the converse property: any optimal plan is actually induced by a map. In this paper, we implement this strategy for a class of metric spaces called Alexandrov spaces.

An Alexandrov space is the natural generalization of Riemannian manifolds whose sectional curvature is bounded below (a precise definition is given in Section 2; we do note however that an Alexandrov space is not necessarily a Gromov-Hausdorff limit of Riemannian manifolds). In their seminal paper \cite{bgp}, Burago, Gromov, and Perelman showed that such (finite dimensional) metric spaces possess a certain kind of "Riemannian structure". This viewpoint was developed further by Otsu and Shioya \cite{otsu}. We use their results to extend McCann's theorem to this setting. Our main result is the following
\begin{theorem}\label{main}
Let $(X,d)$ be a finite dimensional Alexandrov space and $\mu_H$ be the corresponding Hausdorff measure. Let $\mu_0,\mu_1$ be probability measures on $X$ with compact supports such that $\mu_0$ is absolutly continuous with respect to $\mu_H$.

Under these assumptions, Kantorovitch problem defined above admits a solution, and any optimal plan is supported in the graph of a Borel function $F$. This map $F$ is also a minimizer of Monge's problem and satisfies for $\mu$ almost every $x \in X$,  
$$ F(x) = exp (-\nabla \phi(x)), $$
where $\phi$ is a $d^2$-concave function (see Definition \ref{cost}).

Moreover, up to modifications on negligeable sets, the map $\nabla \phi$ is unique, and hence so is the optimal map $F$.

\end{theorem}  

\begin{remark}Observe that in the above theorem, the lower bound on the (Alexandrov) curvature does not  appear explicitly in the statement. Consequently, our result also applies to any compact Riemannian manifold and allows us to give another proof of McCann's theorem \cite{mccann1}. The original proof, relying strongly on the regularity of the Riemannian exponential map, cannot be adapted to this case.
\end{remark}  

Let us mention a geometric motivation which lead us to study the problem above. Recently,  the notion of having bounded below Ricci curvature has been extended to the setting of metric measure spaces independently by Lott-Villani, and by Sturm  \cite{LV1,LV2,St1,St2}. This equivalent definition uses optimal mass transport theory. For a Riemannian manifold, the proof of the equivalence between both definitions basically relies on McCann's theorem \cite{mccann1} and on a non-smooth change of variables formula proved in \cite{mccann4}.  A natural question (mentioned in \cite{LV1}) is whether or not an Alexandrov space has Ricci curvature bounded below in this generalized sense. Our result should be useful on proving this.

To put our result in perspective, we end this section by citing known results about Monge's problem in a singular setting (all but one for  $\mathbb{R}^n$, the other result is for a compact manifold). First, Gangbo and McCann \cite{mccann2} extend Brenier's result to the case of strictly convex cost functions (including the case of the square of an arbitrary norm). Ambrosio and Rigot treated the case of Heisenberg group \cite{AR}. Ambrosio, Kirchheim and Pratelli proved the existence of optimal maps for the original Monge's problem for crystalline norms \cite{AKP}. Recently, Bernard and Buffoni proved the existence of optimal maps on Finsler manifolds with strictly convex norms \cite{bb}. We refer to the books \cite{vil03, RR, vil07} for more on optimal mass transport.

The rest of this paper is organised as follows. In the next section, we present properties of Alexandrov spaces, refering to the book \cite{bbi} for a background of the theory. We then use these results in Section 3 to prove our main theorem. A sketch of our proof is given at the beginning of this section. In the last section, we indicate  how to adapt our result to other strictly convex costs and how to relax the compactness assumption on the supports of the given measures.

\section{Properties of Alexandrov spaces}
In this section, we summarize known properties of Alexandrov spaces used in the rest of the paper. These results are taken from a paper by Burago, Gromov, and Perelman \cite{bgp} and a paper by Otsu and Shioya \cite{otsu}. We refer to the book \cite{bbi} as textbook on metric spaces; most of the results of \cite{bgp} and proofs can also be found there.

\begin{definition}[Alexandrov space]\label{alex} Let $S^2_k$ be the 2-dimensional \\space form of curvature $k$ and $\delta_k$ be the metric induced by the Riemannian metric. A finite dimensional Alexandrov space $X$ of curvature bounded below by $k$ is a complete, connected, locally compact, geodesic space such that  for all geodesics $\gamma$ in $X$  and $\overline{\gamma}$ in 
$S^ 2_k$  such that their lengths are equal, $L(\gamma)=L(\overline{\gamma})$, $d(p,\gamma(0))=\delta_k(\bar{p},\overline{\gamma}(0))$ and $ d(p,\gamma(1))=\delta_k(\bar{p},\overline{\gamma}(1))$,  the inequality below is satisfied for all $t \in [0,1]$: 
$$ d(p,\gamma (t)) \geq \delta_k(\bar{p},\overline{\gamma}(t)).$$
\end{definition}
 
The lower bound on the curvature of an Alexandrov space allows to define the angle between two geodesics starting  at the same point. Using this property, we can prove the first variation formula for distance functions.
\begin{lemma}[First variation formula]\label{fvf} Let $a$ and $x$ be two distinct points of an Alexandrov space $X$ and $\gamma$  be a unitary geodesic starting from $x$. Then, the equality
$$  d(a,\gamma (t)) = d(a,x) -t\cos \angle_{min}+o (t),$$
holds for all nonnegative number $t$ where $ \angle_{min}$ is the smallest angle between $\gamma$ and a geodesic between $x$ and $a$. 
\end{lemma}
For a proof, we refer to \cite[Corollary 4.5.7 and Remark 4.5.12]{bbi} or \cite[Theorem 3.5]{otsu}.


To prove our main result, we will use the fact that an Alexandrov space is not far from being a (Riemannian) manifold. We explain this point of view, starting with a definition.

\begin{definition} Let $X$ be an $n$-dimensional Alexandrov space. A point $p \in X$ is said to be regular if the tangent cone at $p$ is isometric to Euclidean space and singular otherwise. Throughout the rest of this paper, we denote by $Reg(X)$ (respectively $Sing(X)$) the set of regular (respectively singular) points of $X$. 
\end{definition}
\begin{remark}The curvature bound implies that the tangent cone at each point is unique, and is isometric to the Euclidean cone over the space of directions at $x$. The Hausdorff dimension of any tangent cone is equal to the dimension of $X$. (See \cite{bbi} for a proof.)
\end{remark}

We gather together the main properties of the regular set of an Alexandrov space in three theorems. The first one states that, in a certain sense, the regular set covers almost all of $X$. The second theorem shows there is a  certain kind of differential structure on $Reg(X)$. The last establishes the existence of a Riemannan structure on $Reg(X)$ and the compatibility of this structure with the Alexandrov metric. The order in which these theorems are given below is not chronological; furthermore, we have mixed results proved by various authors. We hope that this non-standard presentation helps the readers understanding.


\begin{theorem}\label{reg}The subset $Reg(X)$ is a  dense measurable set (because it is the intersection of countably many dense open sets) of full measure in $X$. More precisely, the Hausdorff dimension of the singular set satisfies $dim_H (Sing(X)) \leq n-1$.
\end{theorem}   
\begin{proof}
The first statement was proved by Burago-Gromov-Perelman in \cite{bgp}, see also \cite[Chapter 10]{bbi}. The Hausdorff measure property was obtained by Burago-Gromov-Perelman \cite{bgp} and Otsu-Shioya \cite[Theorem A]{otsu} independently.
\end{proof} 
 
\begin{theorem}[Charts on $Reg(X)$]\label{diff}
There exists an atlas $(\phi,U_{\phi})_{\phi \in \Phi}$ on $Reg(X)$. In other words, $\forall \phi  \in \Phi$, $\phi: U_{\phi} \mapsto \mathbb{R}^n$ is a map defined on an open set $U_{\phi}$  of $X$ such that
$$ \bigcup_{\phi \in \Phi}U_{\phi} \supset Reg(X).$$
Moreover, if $\phi,\psi \in \Phi$ are such that  $U_{\phi}\cap U_{\psi}\neq \emptyset$, then there exists $V_{\phi} \subset U_{\phi}$, which is dense and of full measure such that 
$\psi \circ \phi^{-1}$ is continously differentiable on $\phi(V_{\phi}\cap V_{\psi}\cap Reg(X))$ (see Definition \ref{vphi} for a precise definition).  
\end{theorem}
\begin{proof}This statement is proved in \cite[Theorem 4.2]{otsu}.
\end{proof}

\begin{remark} Notice however that $Reg(X)$ can be a dense subset of $X$, an example is provided in \cite{otsu}.
\end{remark}  
\begin{theorem}[Riemannian structure on $Reg(X)$]\label{rs}
\mbox{   }\newline
a) There exists a continous Riemannian metric on $Reg(X)$, {\it i.e.} a family $(g_{\phi})_{\phi \in \Phi}$ of maps such that 
$$ g_{\phi} : U_{\phi} \longrightarrow Sym^+(\mathbb{R}^n)$$ 
is a continous map ($ Sym^+(\mathbb{R}^n)$ denotes the set of symmetric positive definite matrices in $\mathbb{R}^n$), and these maps satisfy the usual formula
$$ g_{\phi}=  ^t(d(\phi \circ \psi^{-1}))\,g_{\psi}\,d(\phi \circ \psi^{-1}).$$
b) The Riemannian structure is compatible with the Alexandrov metric, in the following sense: \\
i) For any $x \in Reg(X)$ any chart $\phi \in \Phi$ such that $x \in U_{\phi}$ and any $\delta >0$, there exists a neighbourhood of $x$ in $U_{\phi}$ such that, the map $\phi$ is a bilipschitz homeomorphism with Lipschitz constants smaller than $1+\delta$ on this neighbourhood.\\   
ii) The tangent cone based at a point $x \in Reg(X)$, endowed with the induced metric is isometric to $(\mathbb{R}^n,g_{\phi}(x))$ (assuming that $x \in U_{\phi}$).\\
iii) The metric induced by the Riemannian metric coincides with the original metric.
\end{theorem}
\begin{proof}Statement b)i) was obtained by Burago, Gromov, and Perelman in \cite{bgp}, see also \cite[Chapter 10]{bbi} for a proof. The others were proved by Otsu and Shioya.
\end{proof}

\begin{remark}
Otsu and Shioya also showed that, up to some modifications (taking averages of suitable distance functions), the  natural maps can be made $C^1$ on the whole image  of $U_{\phi}\cap U_{\psi}\cap Reg(X)$ (see \cite[section 5]{otsu} for a statement). However, for our purpose, the previous version is more convenient.   
\end{remark}
%
%
In order to give some idea of the proof of the theorem of Otsu and Shioya, let us give some details of these charts.

\begin{definition}[strained points and natural charts]
Let $(X,d)$ be a n-dimensional Alexandrov space. A point $p$ is called a  $(n,\ep)$ (or simply $\ep$ if the dimension is implicit) strained point if there are $n$ pairs of points  $(x_i,y_i)_{i \in \{1,\cdots ,n\}}$ in $X$ such that $ \forall i,j \, (i\neq j) \;\in \{1,\cdots,n\}$, 
$$ \tilde{\angle}x_ipy_i) > \pi -\ep,$$
$$\tilde{\angle}x_ipx_j) > \frac{\pi}{2} -10\ep,$$
$$\tilde{\angle}x_ipy_j) > \frac{\pi}{2} -10\ep,$$
$$\tilde{\angle}y_ipy_j) > \frac{\pi}{2} -10\ep,$$
where $\tilde{\angle}$ denotes the comparison angle. The collection $(x_i,y_i)$ itself is called a $(n,\ep)$-strainer for $p$.

We denote by $\phi_{x_1,\cdots,x_n}$ (or simply $\phi$ if there is no ambiguity) the following map :
$$ \begin{array}{rcl} \phi  :  U& \longrightarrow &\mathbb{R}^n \\
                             x & \longmapsto &(d(x_1,x),\cdots,d(x_n,x))
\end{array}$$  
\end{definition}

\begin{remark}Strained points were introduced in \cite{bgp}. In particular, any regular point $x$ is an $\ep$-strained point for arbitrary positive $\ep$ (actually, the regular set consists of points which are $\ep$-strained for arbitrary $\ep>0$). It follows that we can consider the above map $\phi$ for any strainer at $x$, and it can be shown that this map induces a chart (in the sense of Theorem \ref{diff}) for $\ep$ sufficiently small. Throughout the rest of this text, we refer to such a map as a ``natural map''.
\end{remark}   

\begin{definition}\label{vphi}
Let $p$ be a point of $X$. We denote by  $V_p$ the set of points $q$  such that there exists a unique geodesic between $p$ and $q$.\\
Let $\phi_{x_1,\cdots,x_n}$ be a natural map. We set 
$$ V_{\phi}=\bigcap_{i=1}^nV_{x_i}. $$
\end{definition}

One of the main ingredients used to prove Theorem \ref{diff} is the following lemma on regularity of distance functions on Alexandrov space, which is of independent interest.

\begin{lemma}[{\cite[Lemma 4.1]{otsu}}\label{lemdif}
]  Let $\phi=\phi_{x_1,\cdots,x_n}$ be a natural map in a neighbourhood of a regular point $p$ and $q$ be an arbitary  point in $X$. The function $d_q \circ \phi^{-1}$ is continously differentiable on $\phi(V_{\phi}\cap V_q)$.

\end{lemma}
An important consequence of the above theorems is that Radema\-cher's theorem on Lipschitz maps holds in this setting.
\begin{corollary}[Rademacher's theorem]\label{rade}
On a $n$-dimensional \\Alexandrov space of curvature bounded below, the usual notion of differentiability (as in the Riemannian setting), (gradient) vector field and first order expansion are well-defined on $Reg(X)$ hence almost everywhere. Moreover, any Lipschitz function is differentiable almost everywhere with respect to the $n$-dimensional Hausdorff measure. 
\end{corollary} 

\begin{proof}
The differentiability statement is a straightforward corollary of the above results. Let us prove Rademacher 's theorem. $Reg(X)$ is a subset of a separable metric space, hence separable. As a consequence, there exists a countable subset $\Phi_{\mathbb{N}}$ of $\Phi$ such that
$$ \bigcup_{\phi \in  \Phi_{\mathbb{N}}} U_{\phi} \supset Reg(X).$$
Thanks to this property and the fact that $Reg(X)$ is of full measure in $X$, it is sufficient to prove that any Lipschitz map $f$ is differentiable almost everywhere on any open set $U_{\phi}$. Now, let $\phi$ be a natural chart defined on $U_{\phi}$. The map $f \circ \phi^{-1}$ is a Lipschitz map (Theorem \ref{rs} b)1), hence differentiable almost everywhere thanks to the usual Rademacher theorem. We conclude the proof noticing that $\phi$ is differentiable on $ V_{\phi}$ (it is a consequence of the first variation formula), a subset of full measure of $X$ (Theorem \ref{diff}).   

\end{proof}

%
%
%
%
%
%

\section{Optimal map on an Alexandrov space}
The goal of this section is to prove Theorem \ref{main}. The proof is divided in several steps. First, we use a compactness argument (namely, Prokhorov's theorem) to prove the existence of an optimal transport plan ({\it i.e.} a solution of Kantorovich problem). Moreover, this plan is related to special maps called $c$-concave functions. This part is known as Kantorovitch duality and holds in the general setting of Polish spaces. The second step, which is the core of the proof, is a proof of the fact that any optimal plan is supported  in the graph of a map $F$ (up to a negligeable set). In the penultimate step, we prove that the map $F$ sends  the initial measure onto the final one. In the last part, we establish the uniqueness of such a map (up to a negligeable set).   


\subsection{Kantorovich duality}
Throughout this paragraph, we refer to $\mu_0,\mu_1$ as Borel probability measures on a complete separable metric space $X$ without any further assumption on their support.

The dual Kantorovich problem is the problem of maximizing the following quantity
$$ J(\phi,\psi) = \int_{X}\phi(x)d\mu_0(x) + \int_X\psi(y)d \mu_1(y)$$
where $\phi$ and $\psi$ are elements of the space $C^b(X)$ of continous bounded functions on $X$ such that $\forall x,y \in X\times X$, 
$$ \phi(x)+\psi(y) \leq c(x,y).$$
Under quite general assumptions, it is possible to show that extrema of both problems coincide, and, moreover, any minimizer of the Kantorovich problem is associated to a pair of maximizers for the dual problem. Before we state this duality theorem, we recall some definitions which will enable us to describe properties of a maximizing pair. 

\begin{definition}\label{cost} A {\it  cost function} on $X$ is a lower semicontinuous function $c \, : U\times V \mapsto \mathbb{R}^+$, where $U,V$ are Borel subsets of $X$. Let $\phi : U \mapsto \mathbb{R}\cup \{+ \infty\}$ be a proper ({\it i.e.} $\neq + \infty$) measurable function. We define $\phi^c : V \mapsto \mathbb{R}\cup \{-\infty\}$,  the $c$-transform of $\phi$ by the formula 
$$  \phi^c(y) = \inf_{x \in U} c(x,y)-\phi(x).$$
We define in a similar manner the $c$-transform of a function defined on $V$ (keeping $c(x,y)$ unchanged because of the possible asymmetry  of $c$) and we call $\phi$ as above, a $c$-concave function if $(\phi^c)^c=\phi$ (in the rest of the paper, we will write $\phi^{cc}$).
\end{definition}

\begin{remark}\label{lipcost} Assume that a cost function is locally Lipschitz ({\it e.g.} \\$c(x,y)=d^2/2(x,y)$) and $U,V$ are compact sets, then any $c$-concave function is Lipschitz also (see \cite{mccann1} for a proof).
\end{remark}

We refer to \cite{vil03} or \cite{RR} for a more detailed analysis of $c$-concavity. There exist numerous versions of the Kantorovich duality. We give here a version borrowed from the book \cite[Theorem 6.1.5]{ags}.
\begin{theorem}[Kantorovich duality]\label{kanto}Let $(X,d)$ be a complete separable metric space, $\mu_0,\mu_1$ be Borel probability measures on $X$ and $c$ be a cost function such that 
$$ \int_{X \times X} c \; d\mu_0 d\mu_1 < +\infty.$$
Then 
\begin{equation}\label{KD}
\sup_{(\psi,\phi) \in C^b(X)\times C^b(X)} J(\psi,\phi)=\min_{\pi \in S(\mu_0,\mu_1)}\int_{X \times X} c(x,y)d\pi(x,y) 
\end{equation}
where $S(\mu_0,\mu_1)$ denotes the set of probability measures whose marginals are $\mu_0,\mu_1$. In addition, there exists a maximiting pair $(\phi,\phi^c) \in L^1(\mu_0)\times L^1( \mu_1)$ for the dual problem and if  $\gamma$ is an optimal plan then
\begin{equation}\label{gam}
 \phi(x)+ \phi^c(y)=c(x,y) \;\;\; \gamma-\mbox{a.e. in } X\times X.
\end{equation}
Conversely, if there exists $\phi \in L^1(\mu_0)$ such that (\ref{gam}) holds, then $\gamma$ is optimal.
\end{theorem}


\subsection{The support of an optimal plan is a graph}\label{graph}
In this section, we prove the main technical result of this paper, namely that the support of an optimal plan is concentrated on the graph of a function. First, we state a lemma on the differentiability of a function along suitable geodesics.

\begin{lemma}[Differentiability along geodesics]\label{TE}Let $f$ be a function on $X$, differentiable at a  point $x \in Reg(X)$ and $\gamma$ be a unitary geodesic defined on $[0,T]$. Then, the Taylor expansion given below holds:
$$ f(\gamma(t)) = f(x) + t \langle \nabla \psi,\gamma'(0)\rangle + o(t).$$
\end{lemma}    

\begin{proof}
Let $\phi=\phi_{p_1, \cdots ,p_n}$ be a natural map around $x$, Theorem \ref{diff} implies that $f \circ \phi^{-1}$ is differentiable in the usual sense. It remains to prove that $\phi \circ \gamma $ is differentiable at $0$. This follows from the first variation formula (Lemma \ref{fvf}) using the fact we can assume (by definition of a natural chart) that $x \in V_{p_i}$.
\end{proof}




\begin{lemma}\label{tang}Let $\psi$ be a $d^2/2$-concave function on an Alexandrov space $(X,d)$. Then
\begin{equation}\label{in1}
 \frac{1}{2} d^2(x,y)\geq \psi(x)+\psi^c(y)
\end{equation}
for all $x,y \in X$. If $x \in Reg(X)$ is a point where $\psi$ is differentiable, then equality holds in (\ref{in1}) if and only if $y=\exp_x(-\nabla \psi (x))$.
\end{lemma}
\begin{remark}\label{re2} $exp_x(-\nabla \psi(x))$ denotes the geodesic whose direction is $\frac{\nabla \psi}{|\nabla \psi|}$ and is parametrized on $[0,|\nabla \psi(x)|]$. In particular,\\
$d(x,\exp_x(-\nabla \psi (x))=|\nabla \psi(x)|$.
\end{remark}
\begin{proof}
The inequality follows from the definition of $\psi^c$. Under the assumptions of Theorem \ref{main}, a compactness argument leads to the existence of a pair which satisfies the equality. So, let us assume that the equality holds in (\ref{in1}) for a pair $(x,y)$. We set $\theta(z)=d^2(y,z)/2$. Let $\gamma(t)$ be a unitary geodesic starting at $x$ and parametrized by $[0,T]$, and assume that $\gamma(T) \in V_x$. By assumption on $x$ and $y$, we have 
\begin{eqnarray*}
 \theta(\gamma(t))-\theta(x) & = & d^2(\gamma(t),y)/2-d^2(x,y)/2\\
                        & \geq &  \psi(\gamma(t)) +\psi(x)\\
                        & \geq & \langle \nabla \psi (x),\gamma'(0)\rangle +  o(t)
\end{eqnarray*}
where we get the first inequality from (\ref{in1}) and the second from Lemma \ref{TE}. Applying the first variation formula (Lemma \ref{fvf}) to $\theta$ yields 
$$ d^2(\gamma(t),y) \leq d^2 (x,y) -2td(x,y)\cos \angle \sigma'(0),\gamma'(0)) +o(t) $$
where $\sigma$ is a geodesic between $x$ and $y$, and  $ \angle \sigma'(0),\gamma'(0))$ is the angle between $\sigma$ and $\gamma$. Letting $t$ go to $0$, we obtain
$$ \langle  \nabla \psi (x)+d(x,y)\sigma'(0), \gamma'(0) \rangle \leq 0.$$
To prove the reverse inequality, one would like to consider geodesic which starting at $x$ and whose direction is $-\sigma'(0)$. Unfortunately, we cannot suppose such a geodesic exists;  however, applying  the above argument to a sequence of geodesics starting at $x$ and whose directions converge to $-\sigma'(0)$ allows to conclude that the function $\theta$ is differentiable at $x$ and
$$ \nabla \psi (x)=-d(x,y)\sigma'(0).$$ 
(We recall that $V_x$ is a dense subset of $X$.) This gives the characterization of $y$ stated in the lemma.                          
\end{proof}
\begin{remark}\label{re1}Note that we also showed that if $x$ is a point such that $\psi$ is differentiable at $x$, there exists a unique geodesic between $x$ and $exp_x(-\nabla \psi(x))$.
\end{remark} 
 Now, we can deduce from the previous lemma that any optimal plan is actually supported in the graph of a Borel function.

\begin{proposition}Let $\Pi_0$ be a minimizer of the variational problem
$$ \min \int_{X\times X} d^2 (x,y) \,d\Pi(x,y)$$ 
among all couplings on $X\times X$ whose marginals are two probability measures $\mu_0,\mu_1$ that satisfy assumptions of Theorem \ref{main}. Then, there exists a measurable function $F$ such that 
$$ \Pi_0 = (Id,F)_{\sharp}\mu_0.$$
Moreover, the map $F$ is defined by the formula

\begin{equation}
 F(x)= exp_x(-\nabla \psi(x))
\end{equation}
and $\psi$ is a $d^2/2$-concave function.
\end{proposition}
\begin{remark}By abuse of language, we will say that $F$ is a minimizer of the above problem.
\end{remark}
\begin{proof} We deduce from Kantorovich duality (Theorem \ref{kanto}) the existence of a pair $(\psi,\psi^c)$ of Lipschitz maps (see Remark \ref{lipcost}) such that
$$ \int_{X \times X} d^2/2(x,y) \,d\Pi_0(x,y) = \int_X \psi(x) d\mu_0(x) + \int_X \psi^c(y) d\mu_1(y)$$
where $(\psi,\psi^c)$ satisfy (\ref{in1}) by definition of the $c$-transform.
Now, Radema\-cher's theorem and Lemma \ref{tang} imply that the support of $\Pi_0$ is concentrated on the graph of the map
$$ F(x) = exp_x (-\nabla \psi(x)).$$
Note that the map $F$ restricted to the subset of $Reg(X)$ of points where $\phi$ is differentiable is a continous map. Hence, the measurability of $Reg(X)$ (Theorem \ref{reg}) entails the measurability of $F$.

It remains to prove that $F_{\sharp} \mu_0=\mu_1$. Let $Z$ be a subset of $Reg(X)$ of full measure, such that $\forall x \in Z, \psi $ is differentiable at $x$. By definition of $\Pi_0$ and by assumption on $\mu_0$, $\Pi_0(Z\times X)= 1$. As a consequence, if $A$ is a Borel set of  $X$, then the following equalities hold 
\begin{eqnarray*}
\mu_1(A) & = & \int_{X \times A} d\Pi_0 = \int_{Z \times A} d\Pi_0 \\
               &  = & \mu_0(Z \cap F^{-1}(A)) = F_{\sharp} \mu_0 (A).
\end{eqnarray*}


\end{proof}

\subsection{Uniqueness property}
It remains to prove the uniqueness of the optimal plan and  optimal map. The results of Subsection \ref{graph} reduce the proof to the case of the optimal map. Therefore, to complete the proof of Theorem \ref{main}, we just need the following result.
\begin{proposition}
Under the assumptions of Theorem \ref{main}, we set  $t(x)=\exp_x(-\nabla \psi(x))$ where $\psi$ is a $d^2/2$-concave function, a solution of Monge's problem. Then, up to modifications on a negligeable subset, the map $\nabla \psi$ is uniquely determined.
\end{proposition}
\begin{proof}
Let $s$ be another solution. Namely, $s$ is a Borel function, mapping $\mu$ to $t_{\sharp}\mu$, such that:
$$ \int_{X}\psi \,d \mu + \int_{X}\psi^c \,dt_{\sharp}\mu = \int_X \frac{1}{2}d^2(x,s(x))\,d\mu.$$ 
As $t_{\sharp}\mu=s_{\sharp}\mu$ by assumption, we get
$$\int_X \psi(x)+\psi^c(s(x))- \frac{1}{2}d^2(x,s(x))\,d\mu(x)=0,$$ 
where the integrand is nonpositive by definition of the $c$-transform, and consequently is equal to $0$ $\mu$ almost everywhere. Lemma \ref{tang} and Rademacher's theorem allow to conclude that $s=t$ almost everywhere. The uniqueness of $\nabla \psi$ follows from Remarks \ref{re2} and \ref{re1}.
\end{proof}

\section{Generalization to other costs}
As in the Riemannian case, our main result can be adapted to other strictly convex costs. Compactness of supports can also be relaxed. The proof of the theorem below is similar to the proof in the Riemannian case, so we only sketch it. We refer to \cite[Theorem 6.2.4]{ags} for a detailed proof in the Euclidean case.

Throughout the section, we consider a strictly convex cost $c(x,y)$ defined by 
$$ c (x,y) = h(d(x,y)),$$
where $h: \mathbb{R^+} \mapsto \mathbb{R^+}$ is a strictly convex and nondecreasing function.

We also use the approximate differential of map which we recall the definition in this setting.
\begin{definition} We say that $f : X \mapsto  \mathbb{R}$ has an approximate differential at $x \in Reg(X)$ if there exists a map $g : X \mapsto  \mathbb{R}$ differentiable at $x$ such that the set $\{f \neq g\}$ has density $0$ at $x$.
\end{definition}

\begin{theorem}
Let $(X,d)$ be a finite dimensional Alexandrov space and $\mu_H$ be the corresponding Hausdorff measure. Let $\mu_0,\mu_1$ be probability measures on $X$ such that  $\mu_0$ is absolutly continuous with respect to $\mu_H$ and
$$ \int_{X \times X} c \; d\mu_0 d\mu_1 < +\infty.$$
Under these assumptions, Kantorovitch problem  admits a solution, and any optimal plan is supported in the graph of a Borel function $F$. This map $F$ is also a minimizer of Monge's problem and satisfies for $\mu$ almost every $x \in X$,\\  
\begin{center}
$ {\displaystyle F(x) = exp (-}\frac{(h')_+^{-1}(|\tilde{\nabla} \phi(x)|)}{|\tilde{\nabla} (\phi(x))|}\displaystyle \tilde{\nabla} \phi(x)), $
\end{center}

if $|\tilde{\nabla} \phi (x)| \neq 0 $ and $F(x)=x$ otherwise. $\phi$ is a $c$-concave function (see Definition \ref{cost}) and $\tilde{\nabla}$ denotes the approximate gradient of $\phi$.

Moreover, up to modifications on negligeable sets, the map $\tilde{\nabla}\phi$ is unique, and as a consequence so is the optimal map $F$.
\end{theorem}  

\begin{proof}
To get the result, we only have to prove an analogue of Lemma \ref{tang}. To this aim, we set $(\phi,\phi^c)$ a maximizing pair of the dual Kantorovich problem. To circumvent difficulties arising from the noncompactness of supports, we use auxilliary maps defined on compact subsets. We fix $ o \in X$, $R$ a positive number, and define 
$$ \phi_R(x)= \inf_{B(o,R)} c(x,y) - \phi^c(y).$$ 
By assumption on $c$, $\phi_R$ is a locally Lipschitz map. Now, the proof of Lemma \ref{tang} gives us the following equivalence, assuming $\phi_R$ is differentiable at $x$:
$$ \phi_R(x) + \phi^c(y) = c(x,y)$$ 
if and only if\\ 
\begin{center}
${\displaystyle y = exp ( -} \frac{(h')_+^{-1}(|\nabla \phi_R(x)|)}{|\nabla (\phi_R(x))|} \displaystyle \nabla \phi_R(x))$\\
\end{center}
where $(h')_+$ denotes the right derivative of $h$ and the equality reads $y=x$ in the case where $|\nabla (\phi_R(x))|=0$. Now, let $\Pi$ be an optimal plan, Kantorovich duality implies
$$ \phi(x) + \phi^c(y) =c(x,y) \;\; \Pi-\mbox{a.e.}$$
Therefore,  for $\mu_0$ almost every $x \in X$, there exists $y$ such that the above equality holds. This entails that $\bigcup_{R \in \mathbb{N}^*} \{\phi=\phi_R\}$ is a subset of full measure. We conclude the proof as in \cite{ags} (note that Lebesgue's theorem on approximate differentiability holds in our setting thanks to the existence of charts on $Reg(X)$ satisfying Property b)i) of Theorem \ref{rs} and Theorem \ref{reg}). 

%
\end{proof}

{ \bf Acknowledgments.} This work was done during a one-year stay at the Scuola Normale Superiore di Pisa. I would like to thank warmly Luigi Ambrosio for useful discussions and remarks and the institute for its hospitality and support. Many thanks to Catriona also for careful rereading.

\bibliographystyle{plain}
\def\cprime{$'$}

\begin{flushleft}
Jérôme Bertrand\\
Institut de Mathématiques\\
Université Paul Sabatier\\
118 route de Narbonne\\
F31062 Cedex 4 Toulouse\\
{\it email: }bertrand@math.ups-tlse.fr
\end{flushleft}

\end{document}